\documentclass[11pt,a4paper]{article}

\usepackage{epic,eepic,epsfig,amsmath,amssymb,amsthm}
\usepackage{t1enc}

 \usepackage[francais]{babel}

\usepackage{color}

\usepackage{bbm}

\def\virgp{\raise 2pt\hbox{,}}

\def\1{\mathbbm{1}}

\def\virgp{\raise 2pt\hbox{,}}
\def\cdotpv{\raise 2pt\hbox{;}}

\def\1{\mathbbm{1}}

\theoremstyle{remark}

\theoremstyle{definition}

\theoremstyle{definition}

\theoremstyle{definition}

\begin{document}

\title{Assessments in Mathematics, undergraduate degree}

\author{Claire David}

\maketitle
\centerline{Universit\'e Pierre et Marie Curie-Paris 6}
\centerline{Laboratoire Jacques-Louis Lions - UMR 7598}
\centerline{Boîte courrier 187, 4 place Jussieu, F-75252 Paris
cedex 05, France}

\vskip 1cm

 \begin{abstract}
In the sequel, we question the validity of multiple choice questionnaires for undergraduate level math courses. Our study is based on courses given in major French universities, to numerous audiences.

 \end{abstract}

\vskip 1cm

 \noindent \textbf{keywords:} multiple choice questionnaires ; evaluation ; math courses ; undergraduate level.

\vskip 1cm

\section*{Introduction}

How to fairly evaluate a large number of undergraduate students in Mathematics? Is the use of multiple choice questionnaires legitimate and desirable ? Are there any other alternatives? \\

Our study, which is intended to be non-exhaustive, first presents an overview of automated assessments, using multiple-choice questionnaires. We also question the ability of such questionnaires to evaluate a student's ability to formulate reasoning.
We then propose an alternative, involving interactive exercises. Finally, we conclude with an example that seemed instructive to us, on a large population a priori reactive to Mathematics.\\

\section{An overview of automated assessments, using multiple-choice questionnaires}

In 2008-2009, 2009-2010, and 2010-2011, at the University Pierre and Marie Curie-Paris 6, in the first year of the bachelor's degree, a mixed assessment, including a first part of writing, reasoning, and a second part, which consisted in a multi-choice questionnaire, was set up in a course of Vectorial calculs, followed by a large number of students: about 800 in semester 1 and 500 in semester 2. This mixed evaluation had been put in place on the one hand, to propose to the students an evaluation corresponding to computational questions, without difficulty or trap, and, on the other hand, to facilitate the correction. MCQs were written using the " alterQCM " package, available under LaTeX \footnote {https://www.ctan.org/tex-archive/macros/latex/contrib/alterqcm}. The treatment was carried out after a scan of the sheets returned by the students.\\

\begin{figure}[h!]
 \center{\psfig{height=7cm,width=9cm,angle=0,file=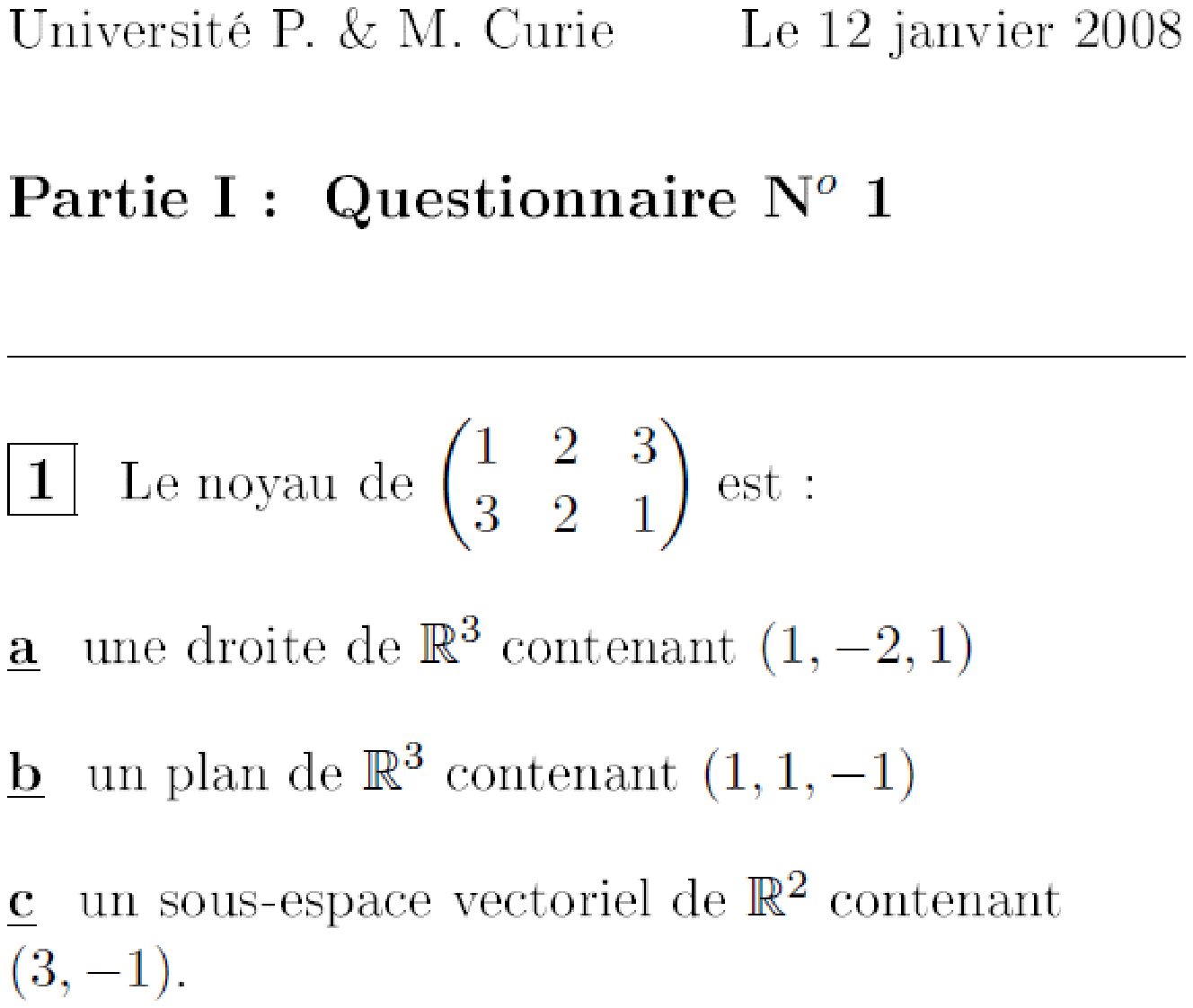}}
 \caption{Extract from the MCQ part of the subject of Vector computing given to UPMC in January 2008.}

\end{figure}

If, compared with more conventional methods, the results were initially better, the teachers quickly realized that this evaluation was not desirable. Indeed, students are often destabilized when going through such evaluations. When they do not what they should answer, they tick a box at random, which has awful consequences due to the counting of negative points in case of wrong answer. \\

We find new occurrences of evaluation by MCQs in the second year of a mathematics degree, in a major French university, in 2014-2015 and 2015-2016.   These assessments reveal many problems. Serious students, with a very correct level otherwise, good or very good, do not pass the exam, having obtained catastrophic results to the MCQs. On the contrary, others achieve surprising results compared to their scores in other courses.

\section{Formulating a reasoning: what is the meaning of automated assessment ? }

It is understandable that the use of automated correction processes is motivated by the high number of students to be assessed, even though inherently the correction of the copies and the writing of the results is part of the teacher's mission.

As is rightly recalled in the study by S.~Bravard (2005), multiple choice questionnaires measure knowledge, understanding. But we can only emphasize the fact that a multiple-choice questionnaire assessment can not assess students' abilities to reason correctly, to write, to argue, to understand.

 Let us give an example of precisely the understanding of limited developments. To carry out a limited development, at a given point, of a given function, consists in approaching this function locally by a polynomial function. But the initial function, and that of the polynomial approximation, are not equal, there is a term of remainder. It is this term of rest on which students often block and make mistakes, errors that are added to those of the actual calculation. A question with pre-programmed answers will not tell if the student really understood this point because it is not him who will write, ONLY, the limited development requested. \\

Let's take advantage of this to highlight the difficulty of setting up a valid MCQ without error. If the MCQ subjects given in the vectorial calculus course were correct,
the MCQ assessments in L2 of Mathematics included errors, with one subject totaling seven questions, three with two correct answers instead of one. On the other hand, some answers were, mathematically, incorrect, which does not correspond to the spirit of a MCQ. This contradicts the basic rules for the elaboration of MCQs which, let us recall, advocate (Leclercq, 1986, Puget, 2010, Missaoui, 2013, Lison, 2014): \\

\begin{enumerate}

\item[$\rightsquigarrow$] To use simple, clear and concise language.

\item[$\rightsquigarrow$] To avoid absolute terms such as "always", "never", "all", "none", or "only ".

\item[$\rightsquigarrow$]   To prefer the positive form.
\end{enumerate}

How, in these conditions, can we evaluate students? \\

The evaluation by MCQ proves to be harmful. Psychologically, already, the students are destabilized by this type of test. Some may argue that the students are no more, or better prepared for a classical subject. This argument does not hold. The French teaching culture is years away from this evaluation format.

In addition, there is a great risk of misleading solutions in the minds of students. Let us refer to Skinner, as S.~Bravard (2005): \\

" Any false solution, in a multiple choice test, increases the probability that a student retrieves a day from his failing memory the incorrect answer instead of the answer
correct. " \\

Moreover, MCQs contribute, as S. Bravard (2005) rightly points out, to contracting the cognitive field of students. \\

Finally, mention may be made of the risks of reading errors of the optical readers. These risks call for manual corrections, not automated ones. The use of multiple-choice questionnaires thus loses some of its usefulness.

 \section{An alternative: Wims}

If we want to use an automated evaluation process, WIMS (Web Interactive Multipurpose Server), created in 1998 by mathematician Xiao Gang \footnote{Unfortunately deceased in 2014.}, a mathematics researcher at the University of Nice, one of the tools currently used most to implement interactive exercises, appears as a hopeful alternative. \\

 Mathematician of genius, but very concerned about teaching, Xiao Gang strongly wanted, beside his research activity, to put in place a pedagogical tool facilitating the assimilation of knowledge, at all levels, from primary to higher education. Wims has been developed under GNU \footnote {Recursive acronym which means " GNU's Not Unix ". GNU is a free operating system project (unlike systems developed by Microsoft), launched by engineering computer scientist Richard Stallman in 1983. From his hacker past, Stallman has kept a fierce will to develop And promote free software.} licence, GPL \footnote{The acronym comes from "General Public License ".}, in OEF \footnote{Open Exercise Format.}. The source code is available, modifiable, distributable. WIMS platforms have been installed locally in many universities \ footnote {Orsay-Paris XI, Paris 1, Paris 6, Paris 7, Paris 12, Paris 13, Caen, Rennes, Marseille, Bordeaux, Grenoble, Lille, University of Savoie, University of the Littoral, ... \\
Source: http://downloadcenter.wimsedu.info/download/map/map.html. Institutions that do not have their own WIMS server choose one of the existing servers.
After the death of Xiao Gang in 2014, WIMS was taken over by Bernadette Perrin-Riou, a professor at University Paris-Sud (Orsay), a specialist in number theory. }. \\

 The resources proposed by WIMS are as follows \footnote{Source: http://wims.unice.fr/description.xhtml}:

\begin{enumerate}
\item[$\rightsquigarrow$] Exercises that can be broken down step-by-step, with the ability to get directions.

\item[$\rightsquigarrow$] Exercises with several possible answers (the answers given are then analyzed by formal calculation software).

\item[$\rightsquigarrow$] Exercises where a counterexample must be given.

\item[$\rightsquigarrow$] Random data drills to work as many times as necessary to do better.
\item[$\rightsquigarrow$] Of the exercises noted \footnote{Wims assigns a note, from 1 to 10.}.
\item[$\rightsquigarrow$] Calculation tools.

\item[$\rightsquigarrow$] Graphics tools.

\item[$\rightsquigarrow$] Course materials, enriched by examples, and related to exercises.

\item[$\rightsquigarrow$] From a more technical point of view, an interface with computing software (Pari, maxima, gap, octave, mupad).
\end{enumerate}

\vskip 0.5cm

From the point of view of the pedagogical organization, WIMS offers the possibility of creating virtual classes allowing \footnote{Source: http://wims.unice.fr/description.xhtml}:

\begin{enumerate}
\item[$\rightsquigarrow$] coaching student work;
\item[$\rightsquigarrow$] taking into account evaluations. The results obtained by the pupils / students are recorded, which allows the teacher to follow the pedagogical progression and, for the learner, to know where it is. From the point of view of analysis of learning, the teacher therefore has an opportunity to trace the actual activity of the learner.
\end{enumerate}

WIMS is constantly evolving, thanks to a community of very active users, who meet regularly to exchange through, among others, the " WIMS Caf\'es ",
where every one can bring his ideas. Originally dedicated to Mathematics, it is now used in many other scientific disciplines: Physics, Chemistry, Life Sciences, ... Interestingly, we note the appearance of exercises in non-scientific disciplines: French (Grammar exercises), Geography.
Interdisciplinary applications are expected. Interestingly, the exercises can be designed either by programming them in their entirety or by using prepared models, allowing teachers unfamiliar with the OEF language to easily make their mark on the WIMS building. From the point of view of mobility, a WIMS Books application for smartphones is under development \ footnote {http://ticewims.unice.fr/mobile/index.html}, as well as another, specific for tablets \ footnote http://ticewims.unice.fr/mobile/index-large.html {}.

\begin{figure}[h!] \center{\psfig{height=5cm,width=10cm,angle=0,file=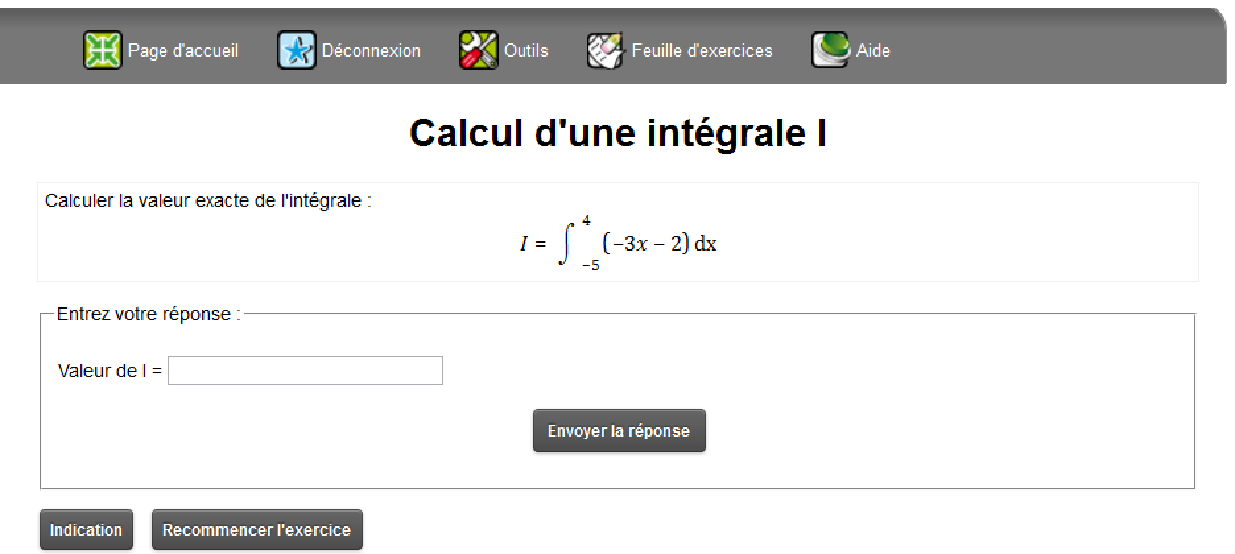}  \psfig{height=5cm,width=10cm,angle=0,file=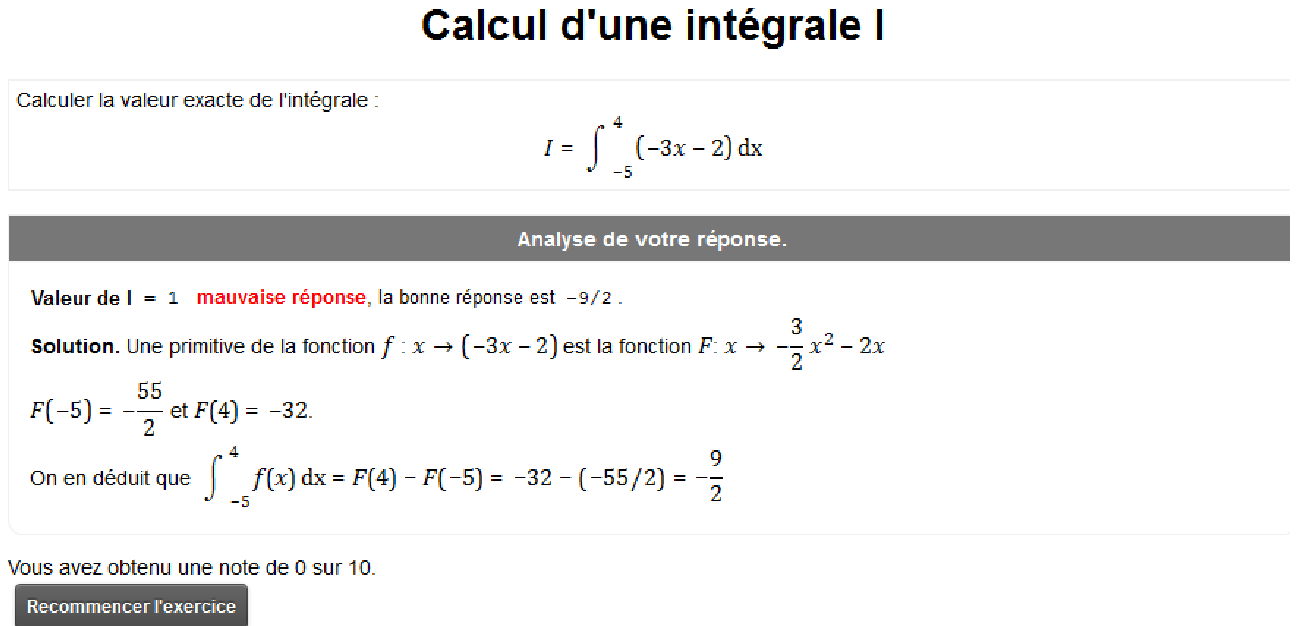}}\\
\centerline{\footnotesize{\textbf{Stated and corrected an interactive exercise of Mathematics Wims (Integral calculus).}}}
\end{figure}

\noindent

\vskip 0.5cm  We begin to have feedback on the contribution of WIMS-type resources, i.e. free access exercise bases, enriched with aids, calculation tools and graphical tools. There is a strong increase in learner motivation, which is reflected in a higher activity (we refer to Hersant (2006) and Cazes (2006)), in line with what was already noted by H.~Ruthven (2002). These studies are corroborated by the many testimonies of teachers working daily with such resources.

 If the evaluation is to be - partially or completely automated, WIMS is therefore a very good alternative. It is a tool (free, in contrast to the MCQs, which represent an expense in terms of correction per optical reader), much more adapted to the offer of training in Mathematics, providing, in addition, a correction enabling the student to Understand its mistakes and progress. Interactivity acts, very positively, in the consolidation of knowledge. WIMS seems to be in perfect harmony with the four pillars of learning as recalled by S.~Dehaene (2012):

\begin{enumerate}
\item[$\rightsquigarrow$] Attention, which is the source of knowledge and action.
\item[$\rightsquigarrow$] Active engagement: the student can not ingest " passively " knowledge.
\item[$\rightsquigarrow$] Feedback: Error signals and individual valuation. The contribution of interactive exercises is essential here. The feedback is similar to that of the learner in gaining feedback. At the same time, taking into account the progress of each one is essential in a process of valorization and encouragement.
\item[$\rightsquigarrow$] The consolidation of knowledge, through the automation of practices, which is achieved through the transfer of the conscious to the unconscious, thus enabling the liberation of resources. Control of basic trigonometry, or limited developments, are exemplary cases of application: after initial learning, regular practice on simple examples allows the student to acquire, in these domains, the " automatisms " required to solve more complex problems: problems of equilibrium or dynamics in point mechanics or solids; Study of asymptotic behavior of functions; Etc ...
\end{enumerate}

\section{A practical study: teaching mathematics in the first year of a bachelor's degree in physics}

\subsection{Context and panorama}

The phenomenon is recent, but speaking: since September 2013, in the undergraduate student population, there is a very notable evolution in terms of the knowledge base, a skill in the field of physics, chemistry, new concepts, ability to reason, logic, misunderstandings. \\

This student public seems to us to be of great pedagogical interest, and for this reason we have chosen to take an in-depth interest in it. Our reference population represents a sample of 519 students. The majority, except~12.9 \%, comes from a general french baccalaur\'eat, branch S (Sciences). \\

Point number one: this population does not like, a priori, mathematics, and chose this path to not make it. Lure, lack of information ? These students do not seem to be aware that the sciences of matter require a minimum of mathematical background, and that if one aspires to pursue at a high level - research in astrophysics, for example, objective of many of them - it is essential to have a very good level of mathematics. \\

Point number two: this population does not seem to have any mathematical training. Among the problems that arise at the beginning of the year are: the resolution of a linear system of two equations to two unknowns; knowledge of the concept of sinus, cosinus; etc ... \\

The task is difficult for teachers. In a semester reduced to twelve weeks (and not $ 6 \times 4 = $ 24 weeks), the students must be provided with the mathematical background giving them the essential knowledge of analysis.

In order to remedy the difficulties mentioned in point one, i.e., a lack of taste, "a priori" for mathematics, and for pedagogical purposes, to re-establish links with physics and chemistry, parts of the course have been arranged in order to synchronize with the teaching of hard sciences: after a first part devoted to the study of numerical functions, we introduce the differential equations, based on numerous examples from Physics: study of electrical circuits, problems of free fall, free oscillations of mass-spring systems, manufacture of soap (saponification of an ester in basic medium, ...). Voluntarily, the emphasis is on the manipulation of the parameters, information jump as described by G.~ Brousseau (1974): \\

 "The information leap consists of finding a fundamental situation that "works" a notion, first choosing the values  of its variables in such a way that the students' previous knowledge makes it possible to develop effective strategies (... ), Then, without changing the rules of the game, change the values of the variables so as to make the complexity of the task to be accomplished much greater. " \\

Thus, the variables are called $ x $, $ t $, $ m $, $ I $, ... and the functions, $ t $, $ x $, $ W $, $ v $ , $ U $, \mbox{$ c $, ...} The task is difficult: if students have encountered many mathematical or physical notions in the course of their previous education, there does not seem to be any more left. As Jean Piaget wrote (Piaget, \& Inhelder, 1968): \\

"Everything is part of the memory (...) outside of which there can be neither understanding of the present nor the same invention. " \\

Alain Lieury (2011) recalls that memorization is necessary for inference reasoning, which contributes to the creation of the meaning ". Semantic memory requires, in order to function optimally, a multiplication of episodes, on a given theme. Recurrence promotes cognitive absorption, and retention of information.

The learning of mathematics can be compared to that of a foreign language. The mastery of the basic tools, as well as that of more sophisticated, is made all the more easily as one regularly practices the discipline. The removal of the misunderstandings that are too often inherent in this subject will be made more easily as the learner is familiar with the handling of the required techniques. From a cognitive point of view, the overly superficial exposure of pupils to multiple notions during schooling has not reached the minimum intensity required, reversing the capacity of absorption.

\subsection{Evaluations}

Faced with heterogeneous audiences, in order to allow a validation of the acquired knowledge accessible to all, and recognizing the aptitudes of each one, the evaluation breaks down into two parts: a part, focused on the skills to be mastered: solving equations, computing differentials, limits, derivatives ; ... not in the form not of a multiple-choice questionnaire, which is a source of errors and is conducive to the " random answers ", but of a answering-document where one is asked to give the solution of basic calculations without justification : in play, only the ability of the student to master a given technique. A second part is more focused on reasoning, and the restoration of formal knowledge. Correct answers to the " skills part ", enable one to pass the exam. The more formal part, more difficult for many, allows a staggering of notes, and a gratification, for those who make the effort. It should be noted that the questions to validate the basic competencies were made, in a slightly different form, if not the same, in progress or in session of tutorials. The required formal knowledge has been listed in advance. The reasoning has been done, in very similar applications. \\

\begin{figure}[h!]
 \center{\psfig{height=8cm,width=11cm,angle=0,file=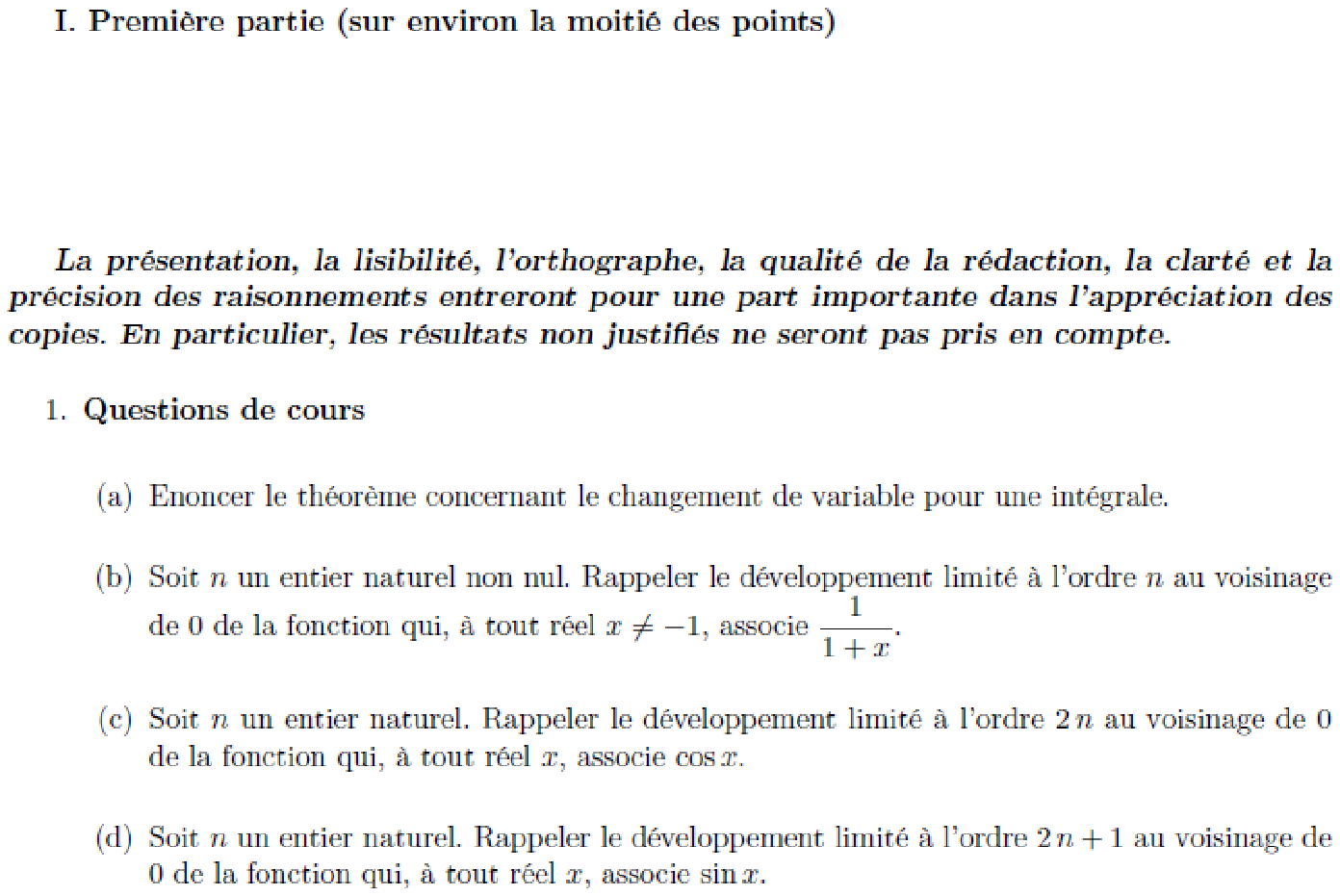}}
 \caption{Extract from the subject of the "Redaction part" of the exam of the Calculus course (UPMC, June 2015).}

\end{figure}

\begin{figure}[h!]
 \center{\psfig{height=8cm,width=11cm,angle=0,file=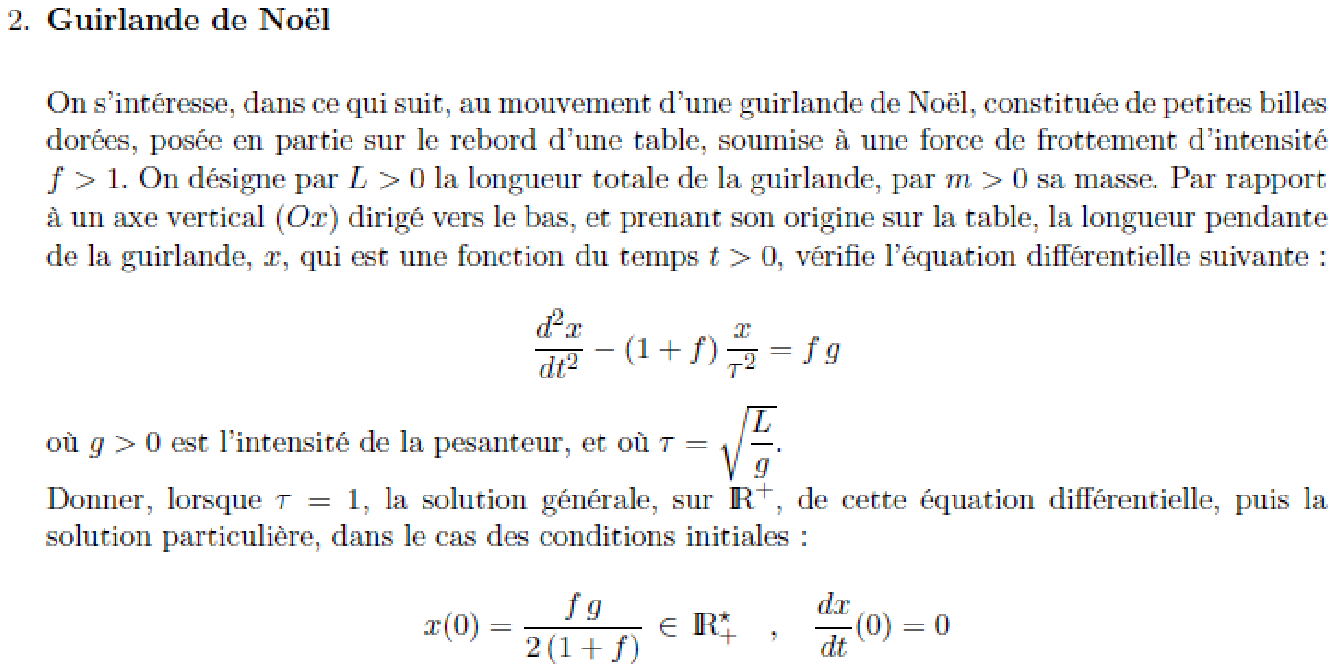}}
 \caption{Extract from the subject of the "Redaction part" of the exam of the Calculus course (UPMC, June 2015).}

\end{figure}

\begin{figure}[h!]
 \center{\psfig{height=7cm,width=11cm,angle=0,file=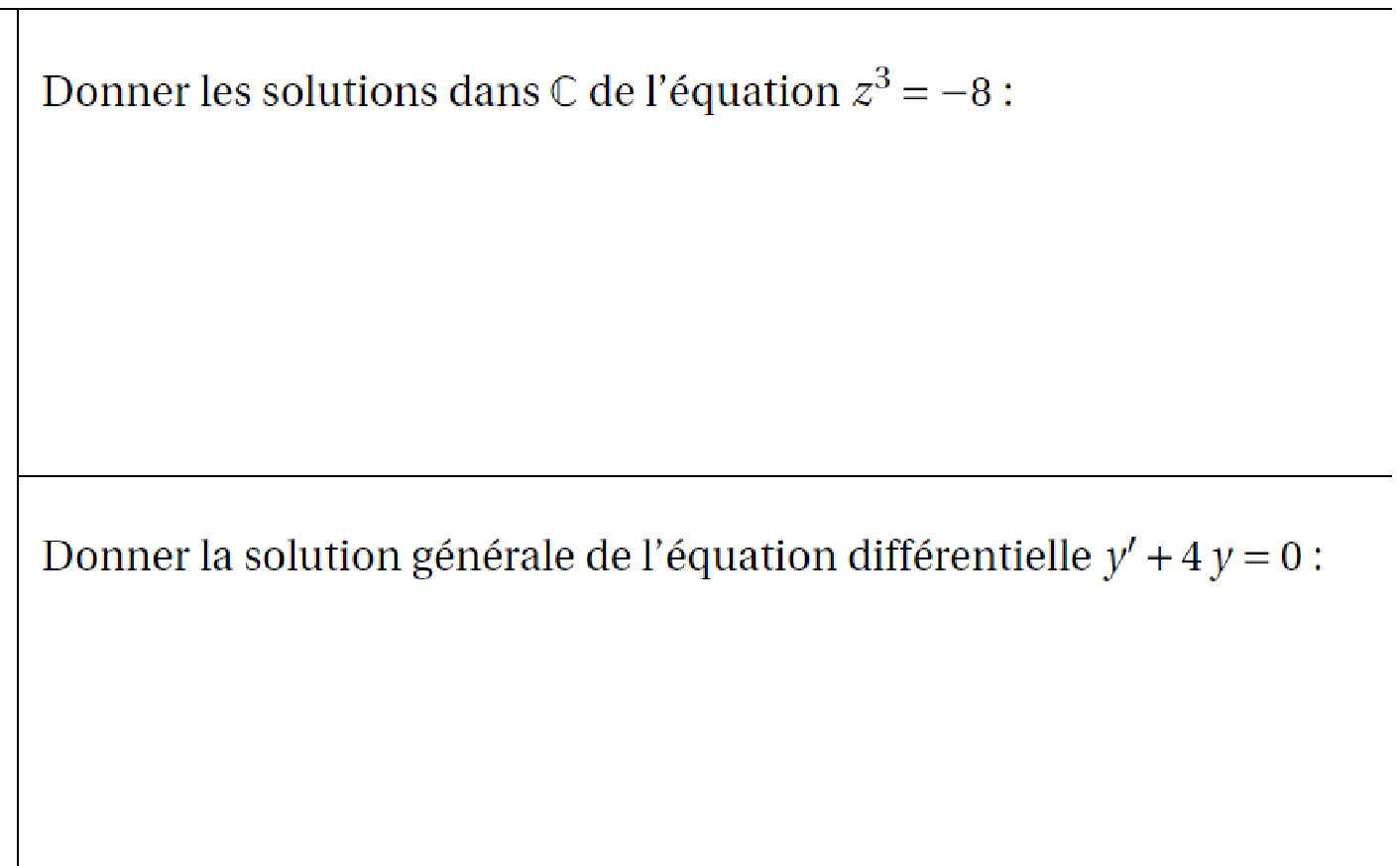}}
 \caption{Extract from the subject of the "Computation part" of the exam of the Calculus course (UPMC, June 2015).}

\end{figure}

These assessments " by competencies ", combined with a weighting, in terms of points, within the course notation, were put in place by the author in 2011-2012, within the framework of vectorial and matrix calculus, in first year of scientific license. The aim was to find a method of assessment that would allow adaptation to a large audience (between 800 and 900 students), very heterogeneous, destined for a second year in physics, chemistry, engineering, electronics or mathematics. From the first iteration (January 2012), this type of evaluation proved to be effective both from the point of view of the teachers and the student public, who, whatever their origin, perceived exactly what was expected of them: Solve a linear system of order two, invert a matrix of order two, compute a determinant, ... \\

We warn against evaluations where the learners of \mbox{" show that ... "}, and where the result is given systematically. Certainly, this prevents a student from moving forward in a subject. All too often, the copies are only a recopying - sometimes inaccurate - of the statements. Students do not always know what to do. In some cases, they embark on a method that has nothing to do with what is expected of them. A mixed evaluation, with a component " per competencies ", breaks the social barrier of an experience that obstructs understanding. Feedback, inducing individual valorisation (Dehaene, 2012), reinforces the learner in his phase of knowledge acquisition, and encourages him to go beyond the simple stage of too basic knowledge. \\

\section*{Conclusion}

Multiple choice questionnaires appear to be totally unsuitable for mathematical evaluations, other than mechanical or purely computational ones, for which, in any case, tools of the type " WIMS " are much better suited. Classical subjects, which evaluate the ability of reasoning, writing, students, with or without a questionnaire with short answers, are, still and always, a hot topic. They alone are the ones most likely to correctly and fairly evaluate students. Of course, human resources are mobilized. This is part of the teaching profession. And, as we know, the scientific entrance examinations at the french "grandes \'ecoles", whether the joint polytechnic competition, the Mines, Centrale, the Ecole Polytechnique and the Ecoles Normales, are still on this format devaluation. There is no secret.

 \end{document}